\documentclass[11pt]{article}

\usepackage{amssymb}

\begin{document}




\begin{center}

{\huge Parabolic partial differential equations with discrete
state-dependent delay: classical solutions and solution manifold}

\bigskip



{{\Large Tibor Krisztin$^{a}$ and Alexander
Rezounenko$^{b,c,}$}\footnote{\footnotesize Corresponding author.
E-mails:  krisztin@math.u-szeged.hu (T. Krisztin),
rezounenko@yahoo.com (A.Rezounenko)} \\

\medskip
$^a$MTA-SZTE Analysis and Stochastic Research Group, Bolyai Institute,\\
University of Szeged,  Aradi v{\'e}rtan{\'u}k tere 1, 6720 Szeged, Hungary \\
$^b$ 
 V.N.~Karazin  Kharkiv National University,  Kharkiv, 61022,
 Ukraine\\
$^c$  Institute of Information Theory and Automation, \\
AS CR, 
P.O. Box 18, 18208 Praha, CR }

\end{center}

\bigskip

\begin{abstract}
Classical solutions to PDEs with discrete state-dependent delay are
studied. We prove the well-posedness in a set $X_F$ which is
analogous to the solution manifold used for ordinary differential
equations with state-dependent delay. We prove that the evolution
operators are $C^1$-smooth on the solution manifold.

\end{abstract}


solution manifold.





\section{ Introduction}
\label{ Introduction}

Differential equations play an important role in describing
mathematical models of many real-world processes. For many years the
models are successfully used to study a number of physical,
biological, chemical, control and other problems. A particular
interest is in differential equations with many variables such as
partial differential equations (PDE) and/or integral differential
equations (IDE) in the case when one of the variables is time. Such
equations are frequently called {\it evolution equations}. They
received much attention from researchers from different fields since
such equations could (in one way or another) discover future states
of a model. It is generally known that taking into account the {\it
past states} of the model, in addition to the present one, makes the
model more realistic. This leads to the so-called delay differential
equations (DDE). Historically, the theory of DDE was first initiated
for the simplest case of ordinary differential equations (ODE) with
constant delay (see the monographs
\cite{Bellman-Cooke_AP-1963_book,Hale,Walther_book,Krisztin-Walther-Wu_AMS-1999_book}
and references therein). Recently many important results have been
extended to the case of delay PDEs with constant delay (see e.g.,
\cite{Travis-Webb_TAMS-1974,Fitzgibbon-JDE-1978,Ruess-1996,Wu_Springer-1996_book}).

Investigating the models described by DDEs it is clear that the
constancy of delays is an extra assumption which significantly
simplifies the study mathematically but is rarely met in the
underlying real-world processes. The value of the delays can be time
or state-dependent. Recent results showed that the theory of
state-dependent delay equations (SDDE) essentially differs from the
ones of constant and time-dependent delays. The basic results on
ODEs with state-dependent delay can be found in
\cite{Driver-AP-1963,Krisztin-Wu-AMUC-1994,Mallet-Paret-Nussbaum-Paraskevopoulos_TMNA-1994,Krisztin-Arino_JDDE-2001,Louihi-Hbid-Arino-JDE-2002,Walther_JDE-2003}
and the review \cite{Hartung-Krisztin-Walther-Wu-2006}. The starting
point of many mathematical studies is the well-posedness of an
initial-value problem for a differential equation. It is directly
connected with the choice of the  space of initial functions (phase
space). For DDEs with constant delay the natural phase space is the
space of continuous functions. However, SDDEs non-uniqueness of
solutions with continuous initial function has been observed in
\cite{Driver-AP-1963} for ODE case. The example in
\cite{Driver-AP-1963} was designed by choosing a non-Lipschitz
initial function $ \varphi \in C[-h,0]$ and a state-dependent delay
such that the value $ - r(\varphi) \in [-h,0]$ (at the initial
function) is a non-Lipschitz point of $\varphi$. In order to
overcome this difficulty, i.e., to guarantee unuique solvability of
initial value problems it was necessary to restrict the set of
initial functions (and solutions) to a set of smoother functions.
This approach includes the restrictions to layers in the space of
Lipschitz functions, $C^1$ functions or the so-called solution
manifold (a subset of $C^1[-h,0]$). As noted in
\cite[p.465]{Hartung-Krisztin-Walther-Wu-2006} 
"...typically, the IVP is uniquely solved for initial and other data
which satisfy suitable Lipschitz conditions." The idea to
investigate ODEs with state-dependent delays in the space of
Lipschitz continuous functions is very fruitful, see e.g
\cite{Mallet-Paret-Nussbaum-Paraskevopoulos_TMNA-1994,Walther_JDE-2003}.
In the present work we rely on the study of solution manifold for
ODEs
\cite{Krisztin-Wu-AMUC-1994,Louihi-Hbid-Arino-JDE-2002,Walther_JDE-2003}

The study of PDEs with state-dependent delay are naturally more
difficult and was initiated only recently
 \cite{Rezounenko-Wu_JCAM-2006,Rezounenko_JMAA-2007,Rezounenko_NA-2009,Rezounenko_NA-2010,Rezounenko_JMAA-2012,Rezounenko-Zagalak-DCDS-A-2013}.
In contrast to the ODEs with state-dependent delays, the possibility
to exploit the space of Lipschitz continuous functions in the case
of PDEs with state-dependent delays meets additional difficulties.
One difficulty is that the solutions of PDEs usually do not belong
to the space of Lipchitz continuous functions. Another difficulty is
that the time-derivative of a solution belongs to a wider space
comparing to the space to which the solution itself belongs. This
fact makes the choice of the appropriate Lipschitz property more
involved, and it depends on a particular model under consideration.
It was already found (see \cite{Rezounenko_NA-2010} and
\cite{Rezounenko-Zagalak-DCDS-A-2013}) that non-local operators
could be very useful in such models and bring additional smoothness
to the solutions. Further studies also show that approaches using
$C^1$-spaces  and solution manifolds (see \cite{Walther_JDE-2003}
and \cite{Hartung-Krisztin-Walther-Wu-2006} for ODE case) could also
be used for PDE models, see
\cite{Rezounenko_NA-2010,Rezounenko-Zagalak-DCDS-A-2013}. In this
work we combine the results for ODEs
\cite{Krisztin-Wu-AMUC-1994,Louihi-Hbid-Arino-JDE-2002,Walther_JDE-2003}
and PDEs \cite{Rezounenko_NA-2010,Rezounenko-Zagalak-DCDS-A-2013}.

We also mention that  a simple and natural additional property
concerning the state-dependent delay  which guarantees the
uniqueness of solutions in the whole space of continuous functions
was proposed in \cite{Rezounenko_NA-2009} and generalized in
 \cite{Rezounenko_JMAA-2012}. We will not develop this approach
 here.

Our goal in this paper is to investigate classical solutions to
parabolic PDEs with discrete state-dependent delay. We find
conditions for the well-posedness and prove the existence of a {\it
solution manifold}. We prove that the evolution operators $G_t :
X_F\to X_F$ are $C^1$-smooth for all $t\ge 0.$ Our considerations
rely on the result \cite{Walther_JDE-2003} and we try to be as close
as possible to the line of the proof in \cite{Walther_JDE-2003} to
clarify which parts of the proof need additional care in the PDE
case. As in \cite{Rezounenko_NA-2010,Rezounenko-Zagalak-DCDS-A-2013}
it is shown that non-local (in space coordinates) operators are
useful in our case. We notice that in
\cite{Rezounenko_NA-2010,Rezounenko-Zagalak-DCDS-A-2013} neither
classical solutions nor $C^1$-smoothness of the evolution operators
were discussed. In the final section we consider an example of a
state-dependent delay which is defined by a threshold condition.

\section{Preliminaries and the well-posedness}
\label{Preliminaries and the well-posedness}

We are interested in the following parabolic partial differential
equation with discrete state-dependent delay (SDD)
\begin{equation}\label{sdd-cl-1}
{d u(t)\over dt} +Au(t) = F\big(u_t \big),\quad t> 0
\end{equation}
with the initial condition
\begin{equation}\label{sdd-cl-ic}
  u_0=u|_{[-h,0]}=\varphi \in C\equiv C([-h,0];L^2(\Omega)). 
\end{equation}
As usual for delay equations \cite{Hale}, for any real $a\le b, t\in
[a,b]$ and any continuous function $u : [a-r,b] \to L^2(\Omega)$, we
denote by $u_t$ the element of $C$ defined by the formula $u_t =
u_t(\theta) \equiv u(t+\theta)$ for $\theta\in [-r,0].$
\smallskip

We assume

{\it {\bf (H1)} Operator $A$ is the infinitesimal generator of a
compact $C_0$-semigroup in $L^2(\Omega).$ }

{\it {\bf (H2)} Nonlinear map $F$ has the form
\begin{equation}\label{F1}
    F(\varphi)\equiv B(\varphi(-r(\varphi))), \qquad F: C\to
    L^2(\Omega),
\end{equation}
where $B: L^2(\Omega)\to L^2(\Omega)$ is a bounded and Lipschitz
operator. Here the state-dependent delay $r: C
([-h,0];L^2(\Omega))\to [0,h]$ is a Lipschitz mapping.
}%

In our study we use the standard (c.f. \cite[def. 2.3,
p.106]{Pazy-1983} and \cite[def. 2.1, p.105]{Pazy-1983})
\smallskip

{\bf Definition~1}. 
 {\it A function $u\in C([-r,T]; L^2(\Omega))$ is called a {\tt mild solution}
 on $[-r,T)$ of the initial value problem (\ref{sdd-cl-1}),(\ref{sdd-cl-ic}) if it satisfies
 (\ref{sdd-cl-ic}) and
 \begin{equation}\label{sdd8-3-1}
u(t)=e^{-A t}\varphi(0) + \int^{t}_0 e^{- A (t-s)}  F(u_s) \, ds,
\quad t\in [0,T).
 \end{equation}
A function $u\in C([-r,T); L^2(\Omega))\bigcap C^1((0,T);
L^2(\Omega))$ is called a {\tt classical solution}  on $[-r,T)$ of
the initial value problem (\ref{sdd-cl-1}),(\ref{sdd-cl-ic}) if it
satisfies  (\ref{sdd-cl-ic}), $u(t)\in D(A)$ for $0<t<T$ and
(\ref{sdd-cl-1}) is satisfied on $(0,T)$.
}%

\medskip

\medskip

{\bf Theorem 1.} {\it Assume (H1)-(H2) are satisfied. Then for any
$\varphi\in C$ there is $t_\varphi>0$
such that initial-value problem (\ref{sdd-cl-1}), (\ref{sdd-cl-ic}) has a 
 mild solution for $t\in [0,t_\varphi)$.}

\smallskip

The proof is standard since $F$ is continuous (see
\cite{Fitzgibbon-JDE-1978}).

We notice that $F$ is not a Lipschitz mapping from $C$ to $[0,h]$,
so we cannot, in general, guarantee the uniqueness of mild solutions
(for ODE case see \cite{Driver-AP-1963}).

\medskip

Let us fix $u$ any mild solution of (\ref{sdd-cl-1}),
(\ref{sdd-cl-ic}) and consider
\begin{equation}\label{sdd-cl-2}
g(t)\equiv  F\big(u_t \big),\quad t\ge 0.
\end{equation}
Mapping $g$ is continuous (from $[0,t_\varphi)$ to $L^2(\Omega)$)
since $B, u$ and $r$ are continuous. Choose $T \in (0,t_\varphi)$.
We have $g\in
C([0,T];L^2(\Omega))$, hence $g\in L^2(0,T;L^2(\Omega))$. The 
initial value problem
\begin{equation}\label{sdd-cl-3}
{d v(t)\over dt} +Av(t) = g(t), \quad v(0)=x\in L^2(\Omega)
\end{equation}
has a unique mild solution, which is $v=u$ if we choose $x=u(0)$.

Now we assume that

\smallskip

{\it {\bf (H3)} operator $A$ is the infinitesimal generator of an
{\bf analytic} (compact) semigroup in $L^2(\Omega).$}

\smallskip

Below we always assume that (H1)-(H3) are satisfied.

As usual, we denote the family of all H\"{o}lder continuous
functions with exponent $\alpha\in (0,1)$ in $I\subset R$ by
$C^\alpha(I;L^2(\Omega))$. By \cite[theorem 3.1, p.110]{Pazy-1983}
the solution $v$ ($=u$) of (\ref{sdd-cl-3}) is H\"{o}lder continuous
with exponent $1/2$ on $[\varepsilon,T]$ for every $\varepsilon \in
(0,T)$. If additionally $x\in D(A)$ then $v\in C^{1\over
2}([0,T];L^2(\Omega))$.

Now we show that $g\in C^{1\over 4}([0,T];L^2(\Omega))$ if
$\varphi\in C^{1\over 2}([-h,0];L^2(\Omega))\subset C$.
Since for $u\in C^{1\over 2}([-h,T];L^2(\Omega))$ and $t\in [0,T]$
one has $||u_t - u_s||_C \le H_u |t-s|^{1\over 2}$ and
$$||g(t)-g(s)|| \le L_B ||u(t-r(u_t)) - u(s-r(u_s))|| \le L_B H_u |t-s + r(u_t) - r(u_s)|^{1\over 2} $$
\begin{equation}\label{sdd-cl-4}
\le L_B H_u \left(|t-s| + L_r ||u_t - u_s||_C \right)^{1\over 2}.
\end{equation}

Here $H_u$ is the H\"{o}lder constant of $u$ on $[-h,T]$, $L_B$ and
$L_r$ are Lipschitz constants.  We get from (\ref{sdd-cl-4}) that
$$||g(t)-g(s)|| \le L_B H_u \left( (T^{1\over 2} +
L_rH_u)|t-s|^{1\over 2} \right)^{1\over 2} \le L_B H_u \left(
T^{1\over 2} + L_rH_u\right)^{1\over 2} |t-s|^{1\over 4}.$$ Here we
used $|t-s|\le T^{1\over 2} |t-s|^{1\over 2}.$ We have shown that
$g\in C^{1\over 4}([0,T];L^2(\Omega))$. It gives, by \cite[corollary
3.3, p.113]{Pazy-1983}, that our mild solution $u$ is {\it
classical} (under assumptions $\varphi\in C^{1\over
2}([-h,0];L^2(\Omega))\subset C$ and $u(0)\in D(A)$).

Set
\begin{equation}\label{sdd-cl-5}
X\equiv \left\{ \varphi\in   C^1 ([-h,0];L^2(\Omega)), \,
\varphi(0)\in D(A) \right\},
\end{equation}
\begin{equation}\label{sdd-cl-8}
||\varphi||_X\equiv \max_{\theta\in [-h,0]}||\varphi (\theta)|| +
\max_{\theta\in [-h,0]}||\dot \varphi (\theta)|| +||A\varphi(0)||.
\end{equation}
Clearly, $X$ is a Banach space since $A$ is closed. We show that
problem (\ref{sdd-cl-1}), (\ref{sdd-cl-ic}) has a {\it unique}
solution for any $\varphi\in X$.

As mentioned before, $F$ is not Lipschitz on $C$, but if $\varphi$
is Lipschitz (with Lipschitz constant $L_\varphi$), then one easily
gets the following estimate (see (\ref{F1}))
\begin{equation}\label{sdd-cl-6}
||F(\varphi)-F(\psi)||\le
L_B||\varphi(-r(\varphi))-\psi(-r(\psi))||$$
$$\le L_B(L_\varphi
|r(\varphi)-r(\psi)| +||\varphi-\psi||_C) 
\le L_B(L_\varphi L_r +1) ||\varphi-\psi||_C.
\end{equation}
Here $L_B$ and $L_r$ are Lipschitz constants of maps $B$ and $r$.

By \cite[theorem 3.5, p.114]{Pazy-1983} (item (ii)), $Au$ and
$du/dt$ are continuous on $[0,T]$, so $u$ is Lipschitz from $[-h,T]$
to $L^2(\Omega)$. This property together with (\ref{sdd-cl-6})
imply the uniqueness of solution to
(\ref{sdd-cl-1}),(\ref{sdd-cl-ic}).

The above proves the following

\medskip

{\bf Theorem 2.} {\it Assume (H1)-(H3) are satisfied. Then for any
$\varphi\in X$ there is $t_\varphi>0$ such that initial value
problem (\ref{sdd-cl-1}), (\ref{sdd-cl-ic}) has a unique
 classical solution for $t\in [0,t_\varphi)$.}

\section{Solution manifold}
\label{Solution manifold}

Let $U\subset $ be an open subset of $X$. We need the following
assumption.

{\it {\bf (S)} The map $F:U\to L^2(\Omega)$ is continuously
differentiable, and for every $\varphi \in U$ the derivative
$DF(\phi)\in L_c(X;L^2(\Omega))$ has an extension $D_eF(\phi)$ which
is an element of the space of bounded linear operators
$L_c(X_0;L^2(\Omega) )$, where $X_0=\{ \varphi\in   C
([-h,0];L^2(\Omega)), \, \varphi(0)\in D(A) \}$ is a Banach space
with the norm $||\varphi||_{X_0}= \max_{\theta\in [-h,0]}||\varphi
(\theta)||  +||A\varphi(0)||$. }


Condition (S) is analogous to that of
\cite[p.467]{Hartung-Krisztin-Walther-Wu-2006}.


\smallskip

 Let us consider the subset
 \begin{equation}\label{sdd-cl-MF}
 X_F= \{ \varphi\in C^1 ([-h,0];L^2(\Omega)), \,
\varphi(0)\in D(A), \dot \varphi (0)+A\varphi(0)=F(\varphi)\}
\end{equation}
of $X$. $X_F$ will be called {\it solution manifold} according to
the terminology of \cite{Walther_JDE-2003}. The equation in
(\ref{sdd-cl-MF}) is understood as equation in $L^2(\Omega)$.  We
have the following analogue to \cite[proposition
1]{Walther_JDE-2003}.

\smallskip

{\bf Lemma~1.} {\it If condition (S) holds and $X_F\neq \emptyset$
then $X_F$ is a $C^1$ submanifold of $X$.
}%

\smallskip

{\it Proof of lemma~1}. 
Consider any $\bar \varphi \in X_F\subset X$ (see (\ref{sdd-cl-MF})
and also (\ref{sdd-cl-5})). Choose $b>0$ so large that
$$|| D_e F(\bar \varphi)||_{L_c(X_0;L^2(\Omega))} <b.
$$
Define $a: [-h,0] \ni s \mapsto s e^{bs} \in R.$ Then
$$ a(0)=0,\quad a^\prime (0)=1,\quad |a(s)| \le {1\over e b}\quad (-h\le s \le 0).
$$
Define the closed subspaces $Y$ and $Z$ of $X$ as follows:
$$Y = \{ a(\cdot )y^0 : y^0\in L^2(\Omega) \} \subset X$$
and
$$ Z = \{ \varphi\in X : \dot \varphi (0)= 0\} \subset X.$$
Clearly $Y\cap Z = \{ 0\}$, and $X= Y \oplus Z$.

We can define the projections
$$ P_Y \phi = a(\cdot)\dot \phi(0), \qquad P_Z \phi = \phi - a(\cdot)\dot
\phi(0).$$ Use $\phi=y+z=P_Y\phi +P_Z\phi$.

We define $$G : X=Y\oplus Z\ni \phi \mapsto \dot \phi(0) +A \phi(0)
-F(\phi) \in L^2(\Omega).$$

Clearly $\phi\in X_F \Longleftrightarrow G(\phi)=0.$ For the bounded
linear map $D_Y G(\bar \varphi) \in { L_c} (Y; L^2(\Omega))$ we have
$$ D_Y G(\bar \varphi)y = \dot y(0)+Ay(0) -DF(\bar \varphi)y = y^0 -
DF (\bar \varphi) a(\cdot)y^0 = y^0 - D_e F (\bar \varphi)
a(\cdot)y^0
$$
since $y=a(\cdot)y^0$ for some $y^0\in L^2(\Omega), \dot y(0)=y^0,
y(0)=0.$

Using the choices of $a$ and $b\in R$ we obtain

$$ ||D_Y G(\bar \varphi)y ||_{L^2(\Omega)} \ge ||y^0||_{L^2(\Omega)}
\left( 1- {||D_e F(\bar \varphi)||\over e b}\right) \ge {1\over 2}
||y^0||_{L^2(\Omega)}.
$$
Then $D_Y G(\bar \varphi):Y\to L^2(\Omega) $ is a linear
isomorphism. The Implicit function theorem can be applied to
complete the proof of lemma.$\blacksquare$


\medskip

For the convenience of the reader we remind some properties of the
semigroup $\{ e^{-At} \}_{t\ge 0}$.

\medskip

{\bf Lemma~2} \cite[theorem 1.4.3, p.26]{Henry-1981} or
\cite[theorem 2.6.13, p.74]{Pazy-1983}. {\it Let $ A$ be a sectorial
operator in the Banach space $Y$ and $Re\, \sigma (A)>\delta>0$.
Then

{\bf (i)} for $\alpha\ge 0$ there exists $C_\alpha<\infty$ such that
\begin{equation}\label{sdd-cl-24}
||A^\alpha e^{-At}||\le C_\alpha t^{-\alpha} e^{-\delta t} \hbox{
for } t> 0;
\end{equation}

{\bf  (ii)} if $0<\alpha\le 1, \, x\in D(A^\alpha)$,
\begin{equation}\label{sdd-cl-25}
||( e^{-At}-I)x||\le {1\over \alpha} C_{1-\alpha} t^{\alpha}
||A^\alpha x|| \hbox{ for } t> 0.
\end{equation}
Also $C_\alpha$ is bounded for $\alpha$ in any compact interval of
$(0,\infty)$ and also bounded as $\alpha\to 0+$.
}%

\smallskip

{\bf Remark 1.} 
{\it It is important to notice that we can 
 write $||(e^{-At}-I)
 A\varphi(0)||\le ||e^{-At}-I||\cdot ||A\varphi(0)||$, but $||e^{-At}-I|| \not\to~0$ as $t\to 0+$
  because  $e^{-At}$ is not a uniformly continuous semigroup
 since $A$ is {\it unbounded} (see \cite[theorem 1.2,
 p.2]{Pazy-1983}).
}

\smallskip

{\bf Remark 2.} {\it
 We also notice that the (linear) mapping $ D(A)\ni \xi \longmapsto (e^{-At}-I)\xi \in
 C^1([0,T];L^2(\Omega))$ is  continuous, while $ L^2(\Omega)\ni \xi \longmapsto (e^{-At}-I)\xi \in
 C^1((0,T];L^2(\Omega))$ is not.
}

\medskip

We need the following 

\medskip

{\bf Lemma 3} . {\it Let $ A$ be a sectorial operator in the Banach
space $Y$ and $f : (0,T) \to Y$ be locally H\"older continuous with
$\int^\rho_0 || f(s)||\, ds~<~\infty $ for some $\rho>0$. For $0\le
t<T$, define
 \begin{equation}
I_T (f)(t)= {\cal F} (t) \equiv \int^t_0 e^{-A(t-s)} f(s)\, ds.
\end{equation}
Then

{\bf (i)} ${\cal F}(\cdot)$ is continuous on $[0,T)$;

{\bf  (ii)} ${\cal F}(\cdot)$ continuously differentiable on
$(0,T)$, with ${\cal F}(t)\in D(A)$ for $0<t<T,$ and $d{\cal
F}(t)/dt +A {\cal F}(t) =f(t)$ on $0<t<T,$ ${\cal F}(t)\to 0$ in $X$
as $t\to 0+$.

{\bf (iii)} If additionally $f:(0,T)\to Y$ satisfies
$$ ||f(t)-f(s)||\le K(s) (t-s)^\gamma \quad \hbox{ for }\quad  0<s<t<T<\infty,$$
where
$K : (0,T)\to {\mathbb R}$ is continuous with $\int^T_0 K(s)\, ds
<\infty$. Then for every $\beta \in [0,\gamma)$ the function
${\cal F} (t)$
 is  continuously
differentiable ${\cal F}:  (0,T) \to Y^\beta\equiv D( A^\beta)$ with
\begin{equation}\label{sdd12-4-11}
    \left\|\frac{d {\cal F} (t)}{dt} \right\|_\beta \le Mt^{-\beta} ||f(t)||+
    M\int^t_0 (t-s)^{\gamma-\beta-1} K(s)\, ds
\end{equation}
for $0<t<T$. Here $M$ is a constant independent of $\gamma,\beta,
f(\cdot)$.

Further, if $\int^h_0 K(s)\, ds = O(h^\delta)$ as $h\to 0+$, for
some $\delta>0$, then $t\to d{\cal F}(t)/dt $ is locally H\"{o}lder
continuous from $(0,T)$ into $Y^\beta$.

{\bf (iv)} If $f :[0,T] \to Y$ is H\"{o}lder continuous (on the
compact $[0,T]$ the local and global H\"{o}lder properties
coincide), then ${\cal F} \in C^1([0,T];Y)$.
}%

\medskip

{\it Proof of lemma 3.} Items (i) and (ii) are proved in \cite[lemma
3.2.1,p.50]{Henry-1981}. Item (iii) is proved in \cite[lemma 3.5.1,
p.70]{Henry-1981}. The proof of (iv) is contained in the proof of
\cite[theorem 3.5, item (ii), p.114]{Pazy-1983}. We briefly outline
the main steps. Using properties (ii) (i.e. $d{\cal F}(t)/dt +A
{\cal F}(t) =f(t)$ on $0<t<T$) and $f\in C([0,T];Y)$ it is enough to
show that $A{\cal F}$ is continuous at $t=0$. We write ${\cal F}(t)
= \int^t_0 e^{-A(t-s)} [f(s)-f(t)]\, ds + \int^t_0 e^{-A(t-s)}
f(t)\, ds = v_1(t) + v_2(t)$. The property $Av_1 \in C^\gamma
([0,T];Y)$ is proved in \cite[lemma 3.4, p.113]{Pazy-1983}. To show
that $Av_2 \in C ([0,T];Y)$ one uses
$$Av_2(t)= \int^t_0 Ae^{-A(t-s)}
f(t)\, ds = \int^t_0 Ae^{-A\tau} f(t)\, d\tau = \int^t_0 \left\{ -
{d\over d\tau} e^{-A\tau} f(t)\right\}\, d\tau $$ $$ =f(0)- e^{-A t}
f(t) = f(0) - e^{-A t} f(0) + e^{-A t} (f(0)- f(t)).$$ Hence
$||Av_2(t)|| \le ||f(0) - e^{-A t} f(0)|| + ||e^{-A t}|| ||f(0)-
f(t)||\le ||f(0) - e^{-A t} f(0)|| + M ||f(0)- f(t)|| \to 0$ as
$t\to 0+$ due to the continuity of $e^{-A t}$ and $f(t)$. It
completes the proof of lemma~3. $\blacksquare$

\medskip


\medskip

To simplify the calculations we assume the following Lipschitz
property holds

 \begin{equation}\label{sdd-cl-9}
\exists\, \alpha\in (0,1), \exists\, L_{B,\alpha}\ge 0 : \forall
u,v\in L^2(\Omega) \Rightarrow ||A^\alpha (B(u)-B(v))|| \le
L_{B,\alpha} ||u-v||.
\end{equation}

{\bf Remark 3.} {\it It is easy to see that (\ref{sdd-cl-9}) implies
similar property with $\alpha=0$ i.e.
\begin{equation}\label{sdd-cl-12}
\exists\, L_{B,0}\ge 0 : \forall u,v\in L^2(\Omega) \Rightarrow
||B(u)-B(v)|| \le L_{B,0} ||u-v||.
\end{equation}
}%

{ \bf Example 1.} Let us consider $B(u)= \int_\Omega f(x-y)
b(u(y))\, dy$ which is a convolution of a function $f\in
H^1(\Omega)$ and composition $b\circ u$ with $b: R\to R$ Lipschitz.
We use the properties of a convolution (see e.g.
\cite[p.104,108]{Brezis2010}) $(f\star g)(x)=\int_\Omega f(x-y)
g(y)\, dy$, namely $||f\star g||_{L^p}\le ||f||_{L^1} ||g||_{L^p}$
for any $f\in L^1$ and $g\in L^p, 1\le p\le \infty$ and also
$D^\beta (f\star g) = (D^\beta f)\star g$, particularly, $\nabla
(f\star g) = (\nabla f)\star g$ (for details see e.g.
\cite[proposition 4.20, p.107]{Brezis2010}).

If we consider Laplace operator with Dirichlet boundary conditions
$A \sim (-\Delta)_D$, then $||A^{1/2}\cdot||$ is equivalent to
$||\cdot ||_{H^1}$, so $||A^{1/2}(B(u)-B(v))||\le
C_1^2||B(u)-B(v)||^2 + C_1^2 ||\nabla(B(u)-B(v))||^2 \le C^2_1
||f||^2_{L^1}||b(u)-b(v)||^2 + C_1^2 ||\nabla
f||^2_{L^1}||b(u)-b(v)||^2$. Using the Lipschitz property of $b$, we
get (\ref{sdd-cl-9}) with $\alpha=1/2$ and $L_{B,\alpha} = C_1 L_b
(||f||^2_{L^1}+ ||\nabla f||^2_{L^1})^{1/2}.$

\medskip

Using (\ref{sdd-cl-9}) and (\ref{F1}) one easily gets the Lipschitz
property for $F$. Namely, for Lipschitz $\psi$ and Lipschitz SDD $r$
$$||A^\alpha (F(\psi)-F(\chi))|| \le ||A^\alpha (B(\psi(-r(\psi)))-B(\chi(-r(\chi))))||
$$ $$\le L_{B,\alpha} ||\psi(-r(\psi)))-\chi(-r(\chi)))||
$$
\begin{equation}\label{sdd-cl-10}
\le L_{B,\alpha} L_\psi L_r ||\psi-\chi||_C + L_{B,\alpha}
||\psi-\chi||_C = L_{F,\alpha}||\psi-\chi||_C, \quad
L_{F,\alpha}=L_{B,\alpha} (L_\psi L_r +1).
\end{equation}

Using (\ref{sdd-cl-12}), similar to (\ref{sdd-cl-10}) one gets
\begin{equation}\label{sdd-cl-13}
||F(\psi)-F(\chi)||\le  L_{F,0}||\psi-\chi||_C, \quad
L_{F,0}=L_{B,0} (L_\psi L_r +1).
\end{equation}

\medskip

We use all notations of \cite{Walther_JDE-2003}, changing $R^n$ for
$L^2(\Omega)$ when necessary. For example,
 we use the notation $E_T$ (see \cite[p.50]{Walther_JDE-2003})
 \begin{equation}\label{sdd-cl-19}
E_T: C^1([-h,0]) \to C^1([-h,T]), \qquad (E_T \varphi)(t)\equiv
\left[ \begin{array}{ll}
\varphi  (t), &\hbox{ for } t\in [-h,0), \\
\varphi(0)+t\dot\varphi(0)  & \hbox{ for } t\in [0,T].
\end{array}
                      \right.
\end{equation}

On the other hand, some notations should be changed. For example,
for any $\psi\in X_F$
 and $r>0$ we set (remind that $||\cdot ||_X$ is not just $C^1$-norm,
 see (\ref{sdd-cl-5}), (\ref{sdd-cl-8}), (\ref{sdd-cl-MF}))
 \begin{equation}\label{sdd-cl-X-psi-r}
X_{\psi,r}\equiv X_F \bigcap \left\{ \psi +
(C^1([-h,0];L^2(\Omega)))_{X,r} \right\} = \left\{ \psi \in X_F :
||\varphi - \psi||_{X} < r \right\}.
\end{equation}

For $T>0$ (to be chosen below), we split a map $x\in C^1([-h,T])
\equiv C^1([-h,T];L^2(\Omega))$ with $x_0=\varphi\in X_F$ given, as
$x=y+\hat \varphi$, where for short  $\hat \varphi  (t)=(E_T
\varphi)(t)$ is defined in (\ref{sdd-cl-19}).

We look for a fixed point of the following map ($\varphi$ is the
parameter)
 \begin{equation}\label{sdd-cl-RTr}
 R_{T r}(\varphi,y)\equiv \left[ \begin{array}{ll}
    e^{-At}\varphi(0) - \varphi(0)- t\dot\varphi(0)+\int^t_0 e^{-A(t-\tau)}F(y_\tau+\hat \varphi_\tau)\, d\tau, & t\in [0,T], \\
    0 \,  & t\in [-h,0), \\
                      \end{array}
                      \right.
\end{equation}
where $R_{T r} : X_{\psi,r}\times
(C^1_0([-h,T];L^2(\Omega)))_\varepsilon \to
C^1_0([-h,T];L^2(\Omega)),$ and $X_{\psi,r}$ defined in
(\ref{sdd-cl-X-psi-r}).

%
\medskip

{\bf Proposition~1}. {\it 
$R_{T r} : X_{\psi,r}\times (C^1_0([-h,T];L^2(\Omega)))
\to C^1_0([-h,T];L^2(\Omega))$.}

\smallskip

 To prove that
the image of $R_{T r}(\varphi,y)=z$ belongs to
 $C^1_0([-h,T];L^2(\Omega)),$  we notice that $y\in
 C^1([-h,T];L^2(\Omega))$ implies $y+\hat \varphi\in
 Lip([-h,T];L^2(\Omega))$, which together with (\ref{sdd-cl-6}) give
 that $F(y_\tau+\hat \varphi_\tau), \tau\in [0,T]$ is Lipschitz, so \cite[lemma 3.2.1,
 p.50]{Henry-1981} can be applied to the integral term in $ R_{T r}$ (see (\ref{sdd-cl-RTr})).
This gives $z\in C^1(0,T; L^2(\Omega))$.

\smallskip

The property $||z(t)||\to 0$ as $t\to 0+$ is simple. The last step
is to show that $||\dot z(t)||\to 0$ as $t\to 0+$. Using \cite[lemma
3.2.1,
 p.50]{Henry-1981} and property $\varphi\in X_F$, we have
$$\dot z(t)= -Ae^{-At}\varphi(0)-\dot \varphi(0) - A\int^t_0 e^{-A(t-\tau)}F(y_\tau+\hat \varphi_\tau)\,
d\tau + F(y_t+\hat \varphi_t)
$$
$$ = -Ae^{-At}\varphi(0) + A\varphi(0) -F(\varphi)- A\int^t_0 e^{-A(t-\tau)}F(y_\tau+\hat \varphi_\tau)\,
d\tau + F(y_t+\hat \varphi_t).
$$
Hence
\begin{equation}\label{sdd-cl-11}
||\dot z(t)||\le ||(e^{-At}-I) A\varphi(0)|| + || F(y_t+\hat
\varphi_t) - F(\varphi)|| + \left\| A\int^t_0
e^{-A(t-\tau)}F(y_\tau+\hat \varphi_\tau)\, d\tau \right\|.
\end{equation}

\smallskip

The first two terms in (\ref{sdd-cl-11}) tend to zero as $t\to 0+$
since $\varphi(0)\in D(A)$, $e^{-At}$ is strongly continuous, $F$ is
continuous and $||y_t+\hat \varphi_t - \varphi||_C\to 0$ as $t\to
0+$. To estimate the last term in (\ref{sdd-cl-11}) we use
(\ref{sdd-cl-10}) for $\psi=0$ and the property $||A^\alpha
e^{-At}||\le C_\alpha t^{-\alpha} e^{-\delta t}, \alpha\ge 0$
(remind that $e^{-At}$ is analytic and see lemma ~2 and
\cite[theorem 1.4.3, p.26]{Henry-1981}, \cite[theorem 2.6.13,
p.74]{Pazy-1983}). So
$$\left\| A\int^t_0
e^{-A(t-\tau)}F(y_\tau+\hat \varphi_\tau)\, d\tau \right\| = \left\|
\int^t_0  A^{1-\alpha} e^{-A(t-\tau)} A^\alpha F(y_\tau+\hat
\varphi_\tau)\, d\tau \right\| $$
 $$\le \int^t_0 C_{1-\alpha}
(t-\tau)^{\alpha-1} e^{-\delta (t-\tau)} L_{B,\alpha} ||y_\tau+\hat
\varphi_\tau||_C\, d\tau $$ $$\le L_{B,\alpha} C_{1-\alpha} \cdot
\max_{s \in [0,T]} ||y_s+\hat \varphi_s||_C \int^t_0
(t-\tau)^{\alpha-1} e^{-\delta (t-\tau)}  \, d\tau \to 0$$ as $t\to
0+$ since the last integral is convergent for $\alpha >0$. It
completes the proof of
Proposition~1.~$\blacksquare$ 

\smallskip

{\bf Remark 4.} {\it It is important in the proof of Proposition~1
to have the property (\ref{sdd-cl-9}) with $\alpha>0$ for the
convergence of the last integral.
}%

\medskip

As in \cite[p.56]{Walther_JDE-2003} we will use local charts of the
manifold $X_F$ and a version of Banach's fixed point theorem with
parameters (see e.g., Proposition 1.1 of Appendix VI in
\cite[p.497]{Walther_book}).

\smallskip

{\bf Remark 5.} {\it More precisely, we look for a fixed point of
$R_{T r}(\varphi,y)$ as a function of $y$ where parameter is
the image of $\varphi$ under a local chart map instead of
$\varphi\in X_{\psi,r}$. The reason is that the parameter should
belong to an open subset of a Banach space, but $X_{\psi,r}$ is not
even linear (it is a subset of the manifold $X_F$).}

We remind that for short we denoted by $\hat \varphi \equiv E_T
\varphi$, where
 $ E_T \varphi$ is defined in (\ref{sdd-cl-19}).


\smallskip

{\bf Proposition~2. }\cite[prop.~2]{Walther_JDE-2003}. {\it For
every $\varepsilon>0$ there exist $T=T(\varepsilon)>0$ and
$r=r(\varepsilon)$ such that for all $\varphi \in \psi +
(C^1([-h,0];L^2(\Omega)))_r $ and all $t\in [0,T]$,
$$\hat \varphi_t \in \psi + (C^1([-h,0];L^2(\Omega)))_\varepsilon
$$
}

The proof is unchanged as in \cite[proposition 2]{Walther_JDE-2003},
so we omit it here.

\smallskip

Let us denote $M_T>0$ a constant satisfying $||e^{-As}||\le M_T$ for
all $s\in [0,T]$. Now we prove an analogue to \cite[proposition
3]{Walther_JDE-2003}.

\smallskip

{\bf Proposition~3}. 
{\it For all $\varphi\in X_{\psi,r}$ and $y,w\in
(C^1_0([-h,T];L^2(\Omega)))_\varepsilon$ one has
\begin{equation}\label{sdd-cl-7}
||R_{T r}(\varphi,y)- R_{T
r}(\varphi,w)||_{C^1([-h,T];L^2(\Omega))}\le L_{R_{T r}}
 ||y-w||_{C^1([-h,T];L^2(\Omega))},
\end{equation}
where we denoted for short the Lipschitz constant
\begin{equation}\label{sdd-cl-7A}
L_{R_{T r}} \equiv T L_{F,0,\varepsilon} (M_T+1)
 + T^\alpha C_{1-\alpha}  M_T\alpha^{-1} L_{F,\alpha,\varepsilon}
\end{equation}
with $L_{F,\alpha,\varepsilon}=L_{B,\alpha} (\varepsilon L_r +1)$
and $L_{F,0,\varepsilon}=L_{B,0} (\varepsilon L_r +1)$
 (c.f. (\ref{sdd-cl-10}), (\ref{sdd-cl-13})).
}%


\smallskip

{\it Proof of proposition 3}. Using  (\ref{sdd-cl-13}),  we have for
all
$||\psi||_{C^1}\le \varepsilon$ 
$$||F(\psi)-F(\chi)||\le  L_{F,0,\varepsilon}||\psi-\chi||_C, \quad
L_{F,0,\varepsilon}=L_{B,0} (\varepsilon L_r +1).$$ Let $z= R_{T
r}(\varphi,y), v= R_{T r}(\varphi,w)$ for $y,w\in
(C^1_0([-h,T];L^2(\Omega)))_\varepsilon$. For all $t\in [0,T]$, one
gets
\begin{equation}\label{sdd-cl-14}
||z(t)-v(t)|| \le || \int^t_0 e^{-A(t-\tau)} (F(y_\tau+\hat
\varphi_\tau)-F(w_\tau+\hat \varphi_\tau))\, d\tau ||\le TM_T
L_{F,0,\varepsilon}||y-w ||_{-h,T}.
\end{equation}

Next 
$||\dot z(t)-\dot v(t)|| \le ||F(y_t+\hat \varphi_t) - F(w_t+\hat
\varphi_t)|| + || A\int^t_0 e^{-A(t-\tau)} (F(y_\tau+\hat
\varphi_\tau)-F(w_\tau+\hat \varphi_\tau))\, d\tau ||\le
L_{F,0,\varepsilon} ||y_t-w_t||_{C} + \int^t_0 || A^{1-\alpha}
e^{-A(t-\tau)}|| ||A^\alpha (F(y_\tau+\hat
\varphi_\tau)-F(w_\tau+\hat \varphi_\tau))||\, d\tau $. To estimate
the first term we write $$||y_t-w_t||_{C} = \max_{s\in [-h,0]}
||\int^{t+s}_0 (\dot y(\tau)-\dot w(\tau))\, d\tau || \le \int^T_0
||\dot y(\tau)-\dot w(\tau)||\, d\tau $$ $$\le T
||y-w||_{C^1([-h,T];L^2(\Omega))}.$$

For the second term,  as in proposition~1, we use the property
$||A^\alpha e^{-At}||\le C_\alpha t^{-\alpha} e^{-\delta t},
\alpha\ge 0$ (see \cite[theorem 1.4.3, p.26]{Henry-1981} or
\cite[theorem 2.6.13, p.74]{Pazy-1983}), the
 Lipschitz property (\ref{sdd-cl-10}) and calculations
 $\int^t_0 (t-\tau)^{\alpha-1}\, d\tau =t^\alpha/\alpha$ to get
$$\int^t_0 || A^{1-\alpha}
e^{-A(t-\tau)}|| ||A^\alpha (F(y_\tau+\hat
\varphi_\tau)-F(w_\tau+\hat \varphi_\tau))||\, d\tau $$
$$\le
 C_{1-\alpha} T^\alpha \alpha^{-1} M_T L_{F,\alpha,\varepsilon}
 ||y-w||_{-h,T}.
$$
Hence
 $$||\dot z(t)-\dot v(t)|| \le \left\{ T L_{F,0,\varepsilon}
 + T^\alpha C_{1-\alpha}  M_T\alpha^{-1} L_{F,\alpha,\varepsilon} \right\}
 ||y-w||_{C^1([-h,T];L^2(\Omega))}.
$$


The last estimate and (\ref{sdd-cl-14}) combined give
(\ref{sdd-cl-7}). $\blacksquare$ 

\medskip

The following statement is an analogue to \cite[proposition 4 and
corollary 1]{Walther_JDE-2003}.

\medskip

{\bf Proposition~4}. 
{\it Let $\delta>0$ there exist $T=T(\delta)>0, r=r(\delta)>0$, such
that for all $\varphi\in X_{\psi,r}$ $( ||\psi-\varphi||_X\le r) $
one has
$$||R_{T r}(\varphi,0)||_{C^1([-h,T];L^2(\Omega))} < \delta.
$$
Moreover, for a positive  $\varepsilon$ there exist $\delta>0$ (and
$T=T(\delta)>0, r=r(\delta)>0$ as above) and  $\lambda\in (0,1),$
such that $R_{T r}$ (defined in (\ref{sdd-cl-RTr})) maps the subset
$X_{\psi,r}\times (C^1_0([-h,T];L^2(\Omega)))_{\varepsilon}$ into
the closed ball $Cl\, (C^1_0([-h,T];L^2(\Omega)))_{\lambda\varepsilon}\subset (C^1_0([-h,T];L^2(\Omega)))_\varepsilon$. 
}%

\medskip

{\it Proof of proposition 4}. Consider $z\equiv R_{T r}(\varphi,0)$.
We write for $t\in [0,T]$
$$z(t)= e^{-At}\varphi(0)-\varphi(0) - t\dot \varphi(0)+ \int^t_0 e^{-A(t-\tau)} F(\hat \varphi_\tau)\, d\tau
$$
$$= (e^{-At} -I)(\varphi(0)-\psi(0)) + (e^{-At} -I)\psi(0) -t\cdot (\dot\varphi(0)-\dot\psi(0)) - t \dot\psi(0)
$$
\begin{equation}\label{sdd-cl-15}
+  \int^t_0 e^{-A(t-\tau)} \left\{F(\hat \varphi_\tau)-F(\hat
\psi_\tau)  \right\}\, d\tau +\int^t_0 e^{-A(t-\tau)} F(\hat
\psi_\tau)\, d\tau.
\end{equation}
We estimiate different parts of (\ref{sdd-cl-15}) in the following
ten steps.

\smallskip

\par 1. Using the property $||(e^{-At} -I)x||\le {1\over \alpha} C_{1-\alpha} t^\alpha ||A^\alpha x||$ (see \cite[thm 1.4.3]{Henry-1981}) one gets
$$||(e^{-At} -I)(\varphi(0)-\psi(0))|| \le C_{1\over 2} t^{1\over 2} ||A^{1\over 2}(\varphi(0)-\psi(0))||
\le \hat C  t^{1\over 2} ||A(\varphi(0)-\psi(0))|| $$ $$\le \hat C
t^{1\over 2} ||\varphi-\psi||_X.
$$

\smallskip

\par 2. $|| t\cdot (\dot\varphi(0)-\dot\psi(0))|| \le t\cdot ||\varphi-\psi||_X.$

\smallskip

\par 3. $|| \int^t_0 e^{-A(t-\tau)} \left\{F(\hat \varphi_\tau)-F(\hat \psi_\tau)  \right\}\, d\tau || \le
M_T t L_{F,0} \max_{\tau\in [0,t]} ||\hat \varphi_\tau-\hat
\psi_\tau||_C \le M_T t L_{F,0} (1+T) ||\varphi-\psi||_X.$

\smallskip

\par 4. $|| \int^t_0 e^{-A(t-\tau)} F(\hat \psi_\tau)\, d\tau || \le M_T t L_{B,0}
 \max_{\tau\in [0,t]} ||\hat \psi_\tau||_C \le M_T t L_{B,0}(1+T) ||\psi||_X.$

\smallskip

Now we proceed to estimate the time derivative of $z(t)$
$$
\dot z(t)= -A e^{-At}\varphi(0)- \dot \varphi(0)+ F(\hat \varphi_t)
-A \int^t_0 e^{-A(t-\tau)} F(\hat \varphi_\tau)\, d\tau
$$
$$=-A e^{-At}\varphi(0) +A \varphi(0)+ F(\varphi) +F(\hat \varphi_t) -A
\int^t_0 e^{-A(t-\tau)} F(\hat \varphi_\tau)\, d\tau
$$
$$= (e^{-At} -I)A(\psi(0)-\varphi(0)) - (e^{-At} -I)A\psi(0)  $$ $$+
[ F(\hat \varphi_t) -F(\hat \psi_t)] + [ F(\hat \psi_t) -F(\psi) ] +
[ F(\psi) - F(\varphi) ] $$
\begin{equation}\label{sdd-cl-16}
- \int^t_0 Ae^{-A(t-\tau)} \{ F(\hat \varphi_\tau) - F(\hat
\psi_\tau) \}\, d\tau - \int^t_0 Ae^{-A(t-\tau)} F(\hat \psi_\tau)\,
d\tau .
\end{equation}

\smallskip

We use the following

\par 5. $|| (e^{-At} -I)A(\psi(0)-\varphi(0))|| \le (M_T+1) ||\varphi-\psi||_X.$

\smallskip

\par 6. $|| F(\hat \varphi_t) -F(\hat \psi_t)|| \le L_{F,0} \max_{\tau\in [0,t]} ||\hat \varphi_\tau-\hat
\psi_\tau||_C\le L_{F,0} (1+T) ||\varphi-\psi||_X.$

\smallskip

\par 7. $|| F(\varphi) -F(\psi)|| \le L_{F,0} ||\varphi-\psi||_X.$

\smallskip

\par 8. $|| F(\hat \psi_t) -F(\psi)|| \to 0$ as $t\to 0+$ since $\hat \psi$ is continuous from $[-h,T]$ to $L^2(\Omega)$.

\smallskip

\par 9. $ || \int^t_0 Ae^{-A(t-\tau)} \{ F(\hat \varphi_\tau) - F(\hat
\psi_\tau) \}\, d\tau || = ||  \int^t_0 A^{1-\alpha} e^{-A(t-\tau)}
A^\alpha \{ F(\hat \varphi_\tau) - F(\hat \psi_\tau) \}\, d\tau ||$
$$\le \int^t_0 C_{1-\alpha} (t-\tau)^{\alpha-1} e^{-\delta (t-\tau)}
L_{F,\alpha} ||\hat \varphi_\tau - \hat\psi_\tau||_C \, d\tau \le
C_{1-\alpha} L_{F,\alpha} D_{\alpha,T} ||\varphi-\psi||_X,$$ where
$D_{\alpha,T}\equiv \int^T_0 (T-\tau)^{\alpha-1} e^{-\delta
(T-\tau)} \, d\tau, \, \alpha>0.$

\smallskip

\par 10. Similar to the previous case ($L_{B,\alpha}$ instead of $L_{F,\alpha}$)
$$ ||  \int^t_0 Ae^{-A(t-\tau)} F(\hat \psi_\tau)\,
d\tau ||\le C_{1-\alpha} L_{B,\alpha} D_{\alpha,T} ||\psi||_X.
$$
Now we can apply estimates 1.-10. (combined) to (\ref{sdd-cl-15}),
(\ref{sdd-cl-16}). It gives the possibility to choose small enough
$T=T(\delta)>0, r=r(\delta)>0$ such that
\begin{equation}\label{sdd-cl-17}
||z||_{C^1([-h,T];L^2(\Omega))} \equiv ||R_{T
r}(\varphi,0)||_{C^1([-h,T];L^2(\Omega))} < \delta .
\end{equation}

\smallskip

{\bf Remark 6.} {\it Small $r$ is used in 5.-7. only. For all the
other terms it is enough (to be small) to have a small $T$.}

 \smallskip

Now we prove the second part of proposition~4. We have
\begin{equation}\label{sdd-cl-18}
||R_{T r}(\varphi,y)||_{C^1([-h,T];L^2(\Omega))} \le ||R_{T
r}(\varphi,y)-R_{T r}(\varphi,0)||_{C^1([-h,T];L^2(\Omega))}
$$ $$+||R_{T r}(\varphi,0)||_{C^1([-h,T];L^2(\Omega))}.
\end{equation}


The first term in (\ref{sdd-cl-18}) is controlled by proposition~3
(see (\ref{sdd-cl-7})), while the second one by (\ref{sdd-cl-17}).

More precisely, we proceed as follows. First choose $\varepsilon>0$,
then choose small $T(\varepsilon)>0$ to have the Lipschitz constant
$L_{R_{T r}} <1$ (see (\ref{sdd-cl-7}), (\ref{sdd-cl-7A})). Next we
set $\delta\equiv {\varepsilon\over 2} (1-L_{R_{T r}})>0$ and the
corresponding $T=T(\delta)\in (0,T(\varepsilon)], r=r(\delta)>0$ as
in the first part of proposition~4, see (\ref{sdd-cl-17}). Finally,
we set $\lambda\equiv  {1\over 2} (1+L_{R_{T r}}) \in (0,1)$. Now
estimates (\ref{sdd-cl-18}), (\ref{sdd-cl-7}) and (\ref{sdd-cl-17})
show that for any $y \in (C^1_0([-h,T];L^2(\Omega)))_{\varepsilon}$
we have
$$||R_{T r}(\varphi,y)||_{C^1([-h,T];L^2(\Omega))} \le L_{R_{T r}} ||y||_{C^1([-h,T];L^2(\Omega))} + \delta
\le L_{R_{T r}} \varepsilon + \delta
$$
$$=  L_{R_{T r}} \varepsilon + {\varepsilon\over 2} (1-L_{R_{T r}}) =
\varepsilon  {1\over 2} (1+L_{R_{T r}}) = \varepsilon \lambda
<\varepsilon.
$$
It completes the proof of proposition~4.$\blacksquare$ 

\medskip

We assume

\medskip

{\it {\bf (H4)} Nonlinear operators $B : L^2(\Omega)\to D(A^\alpha)$
for some $\alpha>0 $ and $r: C([-h,0];L^2(\Omega))\to [0,h]$ are
$C^1$-smooth.
}%

\medskip

{\bf Remark 7.} 
{\it  
Assumption (H4) implies that the restriction \\ $r:
C^1([-h,0];L^2(\Omega))\to [0,h]$ is also $C^1$-smooth. In addition,
it is easy to see that (H4) implies condition (S).
}%

\medskip

{\bf Proposition~5}. 
{\it Assume (H1)-(H4) are satisfied.
Then $R_{T r}$ is $C^1$-smooth.
}

\medskip

The {\it proof of proposition 5} follows the one of
\cite[prop.5]{Walther_JDE-2003}. The main essential difference is
the following. The $C^1$-smoothness of  $B : L^2(\Omega)\to
D(A^\alpha)$ implies the $C^1$-smoothness of  $\widetilde{F} :
X_{\psi,r}\times C^1([-h,0];L^2(\Omega))\to D(A^\alpha)$ defined as
$\widetilde{F}(\varphi,y) \equiv B (\varphi(-r(\varphi+y)) +
y(-r(\varphi+y)))$.

We also use evident additional property  of the $C^1$-smoothness of
the map $ X\ni \varphi \mapsto e^{-At}\varphi(0) \in
C([0,T];L^2(\Omega))$ (remind the definition of $X$ in
(\ref{sdd-cl-5})). Here we use $I_T :  C^1([0,T];L^2(\Omega)) \to
C^1([0,T];L^2(\Omega))$ given by $I_T(y)(t)\equiv \int^t_0
e^{-A(t-\tau)}y(\tau)\, d\tau$ instead of $I_T$ used in
\cite[p.50]{Walther_JDE-2003}. We rely on \cite[lemma
3.2.1,  p.50]{Henry-1981} (see lemma~3, item (iv) above). $\blacksquare$ 


\medskip

As in  \cite[p.56]{Walther_JDE-2003} we are ready to use local
charts of the submanifold $X_F$ and a version of Banach's fixed
point theorem with parameters  (see e.g, \cite[proposition 1.1 of
Appendix VI]{Walther_book}). Namely, propositions 3-5 allow us to
apply the Banach's fixed point theorem to get for any $\varphi\in
X_{\psi,r}$ the unique fixed point $y=y^\varphi\in
(C^1_0([-h,T];L^2(\Omega)))_{\varepsilon}$ of the map $R_{T r}$. We
denote this correspondence by $Y_{T r} : X_{\psi,r} \to
(C^1_0([-h,T];L^2(\Omega)))_{\varepsilon}$ and it is $C^1$-smooth. 

\medskip

It also gives that the map
\begin{equation}\label{sdd-cl-20}
S_{T r} : X_{\psi,r} \to C^1([-h,T];L^2(\Omega)),
\end{equation}
defined by $S_{T r}\varphi  =x^\varphi\equiv y^\varphi+\hat \varphi
\equiv  Y_{T r}(\varphi) + E_T \varphi$  is $C^1$-smooth. Here $E_T
\varphi$ is defined in (\ref{sdd-cl-19}).

\medskip

The local semiflow
$$F_{T r} : [0,T]\times  X_{\psi,r} \to X_F \subset X
$$
is given by
\begin{equation}\label{sdd-cl-21}
F_{T r} (t,\varphi) = x^\varphi_t = ev_t(S_{T r} (\varphi)).
\end{equation}
Here we denoted the evaluation map
\begin{equation}\label{sdd-cl-22}
 ev_t : C^1([-h,T];L^2(\Omega)) \to C^1([-h,0];L^2(\Omega)), \qquad ev_t x \equiv x_t \quad\hbox{for all}\quad t\in [0,T].
\end{equation}

\medskip

{\bf Proposition~6}. 
{\it Assume (H1)-(H4) are satisfied.
Then $F_{T r}$ is continuous, and each solution map $F_{T r}
(t,\cdot):  X_{\psi,r}\backepsilon \phi \mapsto x^{(\phi)}_t \in
X_F, t\in [0,T],$ is $C^1$-smooth. For all $t\in [0,T],$ all $\phi
\in X_{\psi,r}$, and all $\chi \in T_\phi X_F$, one has $T_{F_{T r}
(t,\phi)} \backepsilon D_2 F_{T r} (t,\phi) \chi =
v^{(\phi,\chi)}_t, $  where the function $v\equiv v^{(\phi,\chi)}
\in C^1([-h,T];L^2(\Omega))\cap C([0,T];D(A))$ is the solution of
the initial value problem
\begin{equation}\label{sdd-cl-26}
 \dot v (t) = Av(t) +D F (x^{(\phi)}_t)v_t\quad \mbox{for all} \quad t\in [0,T], v_0=\chi.
\end{equation}

Here $ T_\phi X_F $ is the tangent space to the manifold $X_F$ at
point $\phi\in X_F.$
}

\medskip
{\it Proof of proposition 6}. We denote for short $G\equiv F_{T r}$
and $S\equiv S_{T r}$. Now we discuss the continuity of $F$ (remind
the definition of $X$ in (\ref{sdd-cl-5}) and the norm $||\cdot
||_X$ in (\ref{sdd-cl-8})).
$$||G(s,\chi)-G(t,\varphi)||_X =
||x^\chi_s-x^\varphi_t||_{C^1[-h,0]} + ||A(x^\chi(s)-x^\varphi(t))||
$$
$$\le ||x^\chi_s-x^\varphi_s||_{C^1[-h,0]} +
||x^\varphi_s-x^\varphi_t||_{C^1[-h,0]} +
||A(x^\chi(s)-x^\varphi(s))|| + ||A(x^\varphi(s)-x^\varphi(t))||
$$
\begin{equation}\label{sdd-cl-23}
 \le ||S(\chi)- S(\varphi)||_{C^1[-h,T]} + ||x^\varphi_s-x^\varphi_t||_{C^1[-h,0]} +
||A(x^\chi(s)-x^\varphi(s))|| $$ $$+
||A(x^\varphi(s)-x^\varphi(t))||.
\end{equation}
Consider the third term in (\ref{sdd-cl-23}).
$$||A(x^\chi(s)-x^\varphi(s))|| \le ||e^{-As}A(\chi(0)-\varphi(0))||
$$ $$ + \int^s_0 || e^{-A(s-\tau)} A^{1-\alpha} A^{\alpha}
(F(x^\chi_\tau) - F(x^\varphi_\tau)) ||\, d\tau
$$
$$\le ||\chi-\varphi||_X + C_{1-\alpha}T^\alpha \alpha^{-1} M_T
L_{B,\alpha} (L_{x^\varphi} L_r+1) ||x^\chi-x^\varphi||_{C[-h,T]}
$$
$$ \le ||\chi-\varphi||_X + C_{1-\alpha}T^\alpha \alpha^{-1} M_T
L_{B,\alpha} (L_{x^\varphi} L_r+1) ||S(\chi)-
S(\varphi)||_{C[-h,T]}.
$$
We see that due to the continuity of $S\equiv S_{T r}$ (see
(\ref{sdd-cl-20})) the first and the third terms in
(\ref{sdd-cl-23}) tend to zero when $||\chi-\varphi||_X \to 0.$ The
second term in (\ref{sdd-cl-23}) tends to zero as $|s-t|\to 0$ since
$x\in C^1([-h,T];L^2(\Omega))$. The last term in (\ref{sdd-cl-23})
vanishes due to \cite[Theorem 3.5, item (ii), p.114]{Pazy-1983}
(remind that $x^\varphi(0)\equiv \varphi(0)\in D(A)$). We proved the
continuity of $F$. To verify the differential equation for $v$ (see
(\ref{sdd-cl-26})),
we follow the line of arguments presented in
\cite[p.58]{Walther_JDE-2003}. More precisely, we first verify the
integral equation (\ref{sdd8-3-1}) i.e. show that $v$ is a mild
solution to (\ref{sdd-cl-26}). The only difference in our case is
the presence of the operator $A$ which is linear. Hence it does not
add any difficulties in the differentiability of $S\equiv S_{T r}$
when we define for fixed $\phi \in X_{\psi,r}$, and $\chi \in T_\phi
X_F$ the function $v\equiv DS(\phi)\chi \in
C^1([-h,T];L^2(\Omega)).$ Here $DS$ is understood as the
differential of a map between manifolds (see (\ref{sdd-cl-20}) for
the definition of $S$ and \cite{Arnold-ODE-book-1992} for basic
theory of manifolds). One can see \cite[p.58]{Walther_JDE-2003} that
$v_0=ev_0 DS(\phi)\chi = D (ev_0 \circ S)(\phi)\chi = \chi.$ Here
the evaluation map $ev_t$ is defined in (\ref{sdd-cl-22}). Also for
$t\in [0,T]$ and all $\varphi\in X_{\psi,r}$ one has $ev_t
(S(\varphi)) = ev_t x^{(\varphi)}= x^{(\varphi)}_t = F(t,\varphi),$
which implies (see (\ref{sdd-cl-21}))
$$v_t =ev_t DS(\phi)\chi = D (ev_t \circ S)(\phi)\chi = D_2 F(t,\chi).
$$
To show that $v$ satisfies the integral variant of equation
(\ref{sdd-cl-26}) i.e., it is a mild solution to (\ref{sdd-cl-26}),
we first remind (\ref{sdd-cl-20}) and notation $\hat \varphi (t)=(E_T \varphi)(t)$ 
 (\ref{sdd-cl-19}). It gives for $t>0$
$$ S(\varphi)(t)= x^{(\varphi)}(t) = y^{(\varphi)}+E_T \varphi \equiv  Y_{T r}(\varphi) + E_T \varphi $$
$$= e^{-At}\varphi(0) - \varphi(0)- t\dot{\varphi}(0)+\int^t_0
e^{-A(t-\tau)}F(y_\tau+\hat \varphi_\tau)\, d\tau +
\varphi(0)+t\dot\varphi(0)
$$
$$= e^{-At}\varphi(0)+\int^t_0 e^{-A(t-\tau)}F(y_\tau+\hat \varphi_\tau)\,
d\tau.
$$
Hence
\begin{equation}\label{sdd-cl-28}
S(\phi)(t)= e^{-At}\phi(0)+\int^t_0
e^{-A(t-\tau)}F(x^{(\phi)}_\tau)\, d\tau, \qquad t>0,
\end{equation}
and the definition $v\equiv DS(\phi)\chi \in
C^1([-h,T];L^2(\Omega))$ gives for $t>0$ 

$$ v(t)= (DS(\phi)\chi)(t)= \chi(0)+ \int^t_0 e^{-A(t-\tau)}DF(x^{(\phi)}_\tau)\, v_\tau\,
d\tau.
$$
For more details see \cite[p.58]{Walther_JDE-2003}. So $v$ is a mild
solution to (\ref{sdd-cl-26}).

\smallskip

{\bf Remark 8. } {\it To differentiate the nonlinear term in
(\ref{sdd-cl-28}) we apply the same result  on the smoothness of the
substitution operator as in \cite[p.51]{Walther_JDE-2003}. More
precisely, we consider an open set $U\subset C^1([-h,0];
L^2(\Omega)) $ and the open set
$$U_T \equiv \{\eta\in C([0,T]; C^1([-h,0]; L^2(\Omega))) : \eta (t) \in U \,\hbox{for all}\quad
t\in [0,T]\}.$$ It is proved in \cite[Appendix IV,
p.490]{Walther_book} that the substitution operator $F_T : U_T
\backepsilon \eta \mapsto F\circ \eta \in C([-h,0]; L^2(\Omega)) $
is $C^1$-smooth, with $(DF_T (\eta)\chi)(t)=DF(\eta(t))\chi(t)$ for
all $\eta \in U_T, \chi\in C([0,T]; C^1([-h,0]; L^2(\Omega))), t\in
[0,T].$
}

\medskip

To show that $v$ is classical solution we remind first that
assumption (H4) gives the (local) Lipschitz property
for the Frechet derivative $DF: X \supset U\to L^2(\Omega)$ here
$U\subset X$ is an open set. We remind (see e.g.
\cite[p.466]{Hartung-Krisztin-Walther-Wu-2006}) the form of $DF$
using the restricted evaluation map (not to be confused with the
evaluation map $ev_t$ defined in (\ref{sdd-cl-22}))
$$ Ev : C^1([-h,0];L^2(\Omega)) \times [-h, 0]\backepsilon  (\phi, s)  \mapsto \phi(s) \in L^2(\Omega)$$
which is continuously differentiable, with $D_1Ev(\phi, s)\chi =
Ev(\chi, s)$ and \\ $D_2 Ev(\phi, s)1 = \varphi'(s).$ Hence we write
our delay term $F$ as the composition $F\equiv  B \circ Ev \circ (id
\times (-r))$ (see (\ref{F1})) which is continuously differentiable
from $U$ to $L^2(\Omega)$, with
 $$DF(\phi)\chi = DB(\phi(-r(\phi)))[D_1 Ev(\phi,-r(\phi))\chi - D_2 Ev(\phi,-r(\phi))Dr(\phi)\chi]$$
\begin{equation}\label{sdd-cl-27}
  = DB(\phi(-r(\phi)))[\chi(-r(\phi)) - \phi^\prime (-r(\phi)) Dr(\phi)\chi]
\end{equation}
for $\phi \in U$ and $\chi \in C^1([-h,0];L^2(\Omega)).$

Mappings $B$ and $r$ satisfy (H4) and we remind (see remark 7) 
that our $F$ satisfies the condition similar to (S) in
\cite[p.467]{Hartung-Krisztin-Walther-Wu-2006}.  For an example of a
delay term see below.

The (local) Lipschitz property for the Frechet derivative $DF: X\to
L^2(\Omega)$ and the additional smoothness of the initial function
$\chi \in T_\phi X_F\subset X$ gives the possibility to apply
theorem 2 to show that $v$ is a classical solution to
(\ref{sdd-cl-26}). $\blacksquare$

\medskip

Define the set $\Upsilon = \bigcup_{\phi\in X} [0,t(\phi))\times \{
\phi\} \subset [0,\infty)\times X$ and the map $G: \Upsilon \to X$
given by the formula $G(t,\phi)= x^\phi_t$. 
Propositions 1-6 combined lead to the following

\medskip

{\bf Theorem 3}. 
{\it  Assume (H1)-(H4) are satisfied. Then $G$ is continuous, and
for every $t\ge 0$ such that $\Upsilon_t\neq \emptyset$ the map
$G_t$ is $C^1$-smooth. For every $(t,\phi)\in \Upsilon$ and for all
$\chi\in T_\phi X$, one has $DG_t (\phi)\chi = v_t$ with
$v:[-h,t(\phi))\to L^2(\Omega)$ is $C^1$-smooth and satisfies $\dot
v(t)= Av(t)+D F (G(t,\phi))v_t, $ for $t\in [0,t(\phi)), \,
v_0=\chi$.
}

\section{Example of a state-dependent delay}
 \label{Example of a state-dependent delay}

Consider the following example of the delay term used, for example,
in population dynamics  \cite[p.191]{Kuang-1993-book}. It is the
so-called, threshold condition.

The state-dependent delay $r: C ([-h,0];L^2(\Omega))\to [0,h]$ is
given implicitly by the following equation
\begin{equation}\label{sdd-cl-r1}
    R(r;\varphi)=1,
\end{equation}
where
\begin{equation}\label{sdd-cl-r2}
    R(r;\varphi)\equiv \int^0_{-r} \left(  \frac{C_1}{C_2+\int_\Omega \varphi^2(s,x)\, dx} +C_3 \right)\, ds, \quad C_i>0.
\end{equation}
Since $$D_r R(r(\varphi);\varphi)\cdot Dr(\varphi)\psi + D_\varphi
R(r(\varphi);\varphi) \psi =0$$ and
$$D_r R(r(\varphi);\varphi)\cdot 1 = \left(  \frac{C_1}{C_2+\int_\Omega \varphi^2(-r,x)\, dx} +C_3
\right)\cdot 1 \neq 0, \quad C_i>0,$$
$$ D_\varphi R(r(\varphi);\varphi) \psi =
- \int^0_{-r} \left\{  \frac{C_1}{[C_2+\int_\Omega \varphi^2(s,x)\,
dx]^2} \cdot 2 \cdot \int_\Omega \varphi(s,x) \psi(s,x)\, dx
\right\}\, ds,
$$
we have
$$    Dr(\varphi) \psi = \left(  \frac{C_1}{C_2+\int_\Omega
\varphi^2(-r,x)\, dx} +C_3 \right)^{-1}
$$
\begin{equation}\label{sdd-cl-r3}
\times 
\int^0_{-r(\varphi)} \left\{ \frac{C_1}{[C_2+\int_\Omega
\varphi^2(s,x)\, dx]^2} \cdot 2 \cdot \int_\Omega \varphi(s,x)
\psi(s,x)\, dx \right\}\, ds.
\end{equation}
Now, we substitute the above form of $Dr(\varphi) \psi $ into
(\ref{sdd-cl-27}) and  arrive to

$$ DF(\varphi) \psi =DB(\varphi(-r(\varphi)))\left[ \mbox{} \psi(-r(\varphi)) - \varphi^\prime (-r(\varphi)) \times \right.$$

  $$  \left(  \frac{C_1}{C_2+\int_\Omega \varphi^2(-r,x)\, dx} +C_3
\right)^{-1} $$
\begin{equation}\label{sdd-cl-F2}\times \left.
\int^0_{-r(\varphi)} \left\{ \frac{C_1}{[C_2+\int_\Omega
\varphi^2(s,x)\, dx]^2} \cdot 2 \cdot \int_\Omega \varphi(s,x)
\psi(s,x)\, dx \right\}\, ds \right].
\end{equation}

We see that mapping $r$ satisfies (H4). We also remind (see remark 7) 
that in this example $F$ satisfies the condition similar to (S) in
\cite[p.467]{Hartung-Krisztin-Walther-Wu-2006}, provided operator $B
: L^2(\Omega)\to D(A^\alpha)$ (for some $\alpha>0$) is $C^1$-smooth.

\medskip

 {\bf Acknowledgments.} This work was supported in part by GA
CR under project P103/12/2431.





\bigskip \bigskip

\hfill Version: July 6, 2014

\end{document}